\documentclass[12pt,reqno]{amsart}
\usepackage{amsmath,amsfonts,amssymb,amscd,amsthm,amsbsy,epsf}
\usepackage{graphicx,psfrag,comment}
\numberwithin{equation}{section}
\newtheorem{theorem}{Theorem}

\newtheorem{algorithm}[theorem]{Algorithm}
\theoremstyle{remark}

\theoremstyle{definition}
\newtheorem{defn}{Definition}

\def\Id{\operatorname{Id}}

\def\E{{\mathcal E}}
\def\nat{{\mathbb N}}
\def\real{{\mathbb R}}
\def\torus{{\mathbb T}}
\def\zed{{\mathbb Z}}

\def\M{{\mathcal M}}
\def\G{{\mathfrak G}}

\begin{document}
\title[Fast numerical algorithms]
{Fast iteration of cocyles over rotations and Computation of hyperbolic bundles}
\author{Gemma Huguet,  Rafael de la Llave \and Yannick Sire}
\address{Department of Mathematics, The University of Texas at Austin,
Austin, TX 78712} \email{llave@math.utexas.edu}
\address{Centre de Recerca Matem\`atica, Apartat 50, 08193 Bellaterra (Barcelona), Spain}
\email{gemma.huguet@upc.edu}
\address{Universit\'e Paul C\'ezanne, Laboratoire LATP UMR 6632, Marseille, France}
\email{sire@cmi.univ-mrs.fr}
\begin{abstract}
In this paper, we develop numerical algorithms that use small
requirements of storage and operations for the computation of
hyperbolic cocycles over a rotation. We present fast algorithms for the iteration of the
quasi-periodic cocycles and the computation of the invariant bundles, which is
a preliminary step for the computation of invariant
whiskered tori. 
\end{abstract}



\maketitle
\baselineskip=18pt              
\tableofcontents

\section{Introduction}

The goal of this paper is to describe efficient algorithms to
compute quasi-periodic cocycles over rotations.  We present fast algorithms for the
iteration of cocycles over rotations and for the calculation of
their invariant bundles. The main idea is to use a renormalization
algorithm which allows to pass from a cocycle to a longer cocycle.

The calculation of invariant bundles for cocycles is a preliminary
step for the calculation of whiskered invariant
tori. Indeed, these algorithms require the computation of the
projections over the linear subspaces of the linear cocycle.

\section{Some standard definitions on cocycles}

Given a matrix-valued function \hfill
$M :\torus^\ell \to GL(2d,\real)$ and a vector
$\omega\in\real^\ell$, we define the \emph{cocycle} over the
rotation $T_\omega$ associated to the matrix $M$ by a function $\M :\zed\times\torus^\ell \to GL (2d,\real)$ given by
\begin{equation}
\M (n,\theta) =
\begin{cases}
M (\theta  + (n-1)\omega) \cdots M(\theta)&n\ge 1,\\
\Id &n=0,\\
M^{-1} (\theta + (n+1)\omega) \cdots M^{-1}(\theta)&n\le 1.
\end{cases}
\label{cocycle1}
\end{equation}

Equivalently, a cocycle is defined by the relation
\begin{equation}
\begin{split}
\M(0,\theta) & = \Id,\\
\M(1,\theta) & = M(\theta),\\
\M(n+m,\theta) & = \M(n,T_\omega^m (\theta)) \M(m,\theta).
\end{split}
\label{cocycle2}
\end{equation}

We will say that $M$ is the generator of $\M$. Note that if
$M(\torus^\ell)\subset G$ where $G\subset GL (2d,\real)$ is a group,
then $\M(\zed,\torus^\ell) \subset G$.

The main example of a cocycle in this paper is
$$M(\theta) = (DF\circ K)(\theta),$$
for $K$ a parameterization of an invariant torus satisfying the invariance equation
$$F \circ K= K \circ T_\omega. $$

Other  examples appear in discrete Schr\"odinger
equations \cite{Puig02}. In the above mentioned examples, the
cocycles lie in the symplectic group and in the unitary group,
respectively.

Similarly, given a matrix valued function $M(\theta)$, a
\emph{continuous in time cocycle} $\M(t,\theta)$ is defined to be the
unique solution of
\begin{equation}
\begin{split}
\frac{d}{dt} \M(t,\theta) & = M(\theta + \omega t) \M(t,\theta),\\
\M(0,\theta) & = \Id .
\end{split}
\label{cocycle2cont}
\end{equation}

{F}rom the uniqueness part of Cauchy-Lipschitz theorem, we have the following property
\begin{equation}
\begin{split}
\M(\theta, t+s) & = \M(\theta + \omega t,s) \M(\theta,t),\\
\M(\theta ,0) & = \Id .
\end{split}
\label{cocycle1cont}
\end{equation}

Note that  \eqref{cocycle2cont} and \eqref{cocycle1cont} are the exact
analogues of \eqref{cocycle1} and \eqref{cocycle2} in a continuous context. Moreover, if $M(\torus^\ell)\subset \G$, where $\G$ is
 a sub-algebra of the Lie algebra of the Lie group $G$, then $\M(\real,\torus^\ell) \subset G$.

The main example for us of a continuous in time cocycle will be
$$M(\theta) = (DX\circ K)(\theta),$$
where $K$ is a solution of the invariance equation
$$\partial_\omega K = X \circ K$$
and $X$ is a Hamiltonian vector field. In this
case, the cocycle $\M(\theta,t)$ is symplectic.


\section{Hyperbolicity of cocycles}

One of the most crucial property of cocycles is hyperbolicity (or
spectral dichotomies) as described in \cite{MeyerS89,SackerS74, SackerS76a,
SackerS76b, Sacker78}.

\begin{defn}\label{def:dichotomy}
Given $0<\lambda <\mu$ we say that a cocycle $\M(n,\theta)$ (resp.
$\M(t,\theta)$) has a $\lambda,\mu-$ dichotomy if for every $\theta \in \torus^\ell$ there exist a
constant $c>0$ and a splitting depending on $\theta$,
$$T\real^{2d} = \E^s \oplus \E^u$$
which is characterized by:
\begin{equation}
\begin{split}
&(x_\theta, v)\in \E^s \Leftrightarrow |\M(n,\theta)v|
\le c\lambda^n |v|\ ,\qquad \forall n\ge 0\\
&(x_\theta,v)\in \E^u \Leftrightarrow |\M(n,\theta)v| \le c\mu^{n} |v|\
,\qquad \forall n\le 0
\end{split}
\label{characterization}
\end{equation}
or, in the continuous time case
\begin{equation}
\begin{split}
&(x_\theta,v)\in \E^s \Leftrightarrow |\M(t,\theta)v| \le c\lambda^t |v|\ ,\qquad \forall t\ge 0\\
&(x_\theta,v)\in \E^u \Leftrightarrow |\M(t,\theta)v | \le c\mu^{t}|v|\
,\qquad \forall t\le 0.
\end{split}
\label{characterizationcont}
\end{equation}
\end{defn}

The notation $\E^s$ and $\E^{u}$ is meant to suggest that an
important case is the splitting between stable and unstable
bundles. This is the case when $\lambda <1<\mu$ and the cocycle is said to be hyperbolic. Nevertheless, the theory developed
in this section assumes only the existence of a spectral gap.

In the application to  the computation
of tori,  $M(\theta)=(DF\circ K)(\theta)$ and  $x_\theta = K(\theta)$.
The existence of the spectral gap
 means that at every point of the invariant torus $K(\theta)$
one has a splitting so that the vectors grow with
appropriate rates $\lambda,\mu$ under iteration of
the cocycle. In the case of the invariant torus, it
can be seen that the cocycle is just the fundamental
matrix of the variational equations so that the cocycle
describes the growth of infinitesimal perturbations.

It is well known that the mappings $\theta \to \E_{x_\theta}^{s,u}$ are
$C^r$ if $M(.,)\in C^r$ for $r\in\nat\cup\{\infty,\omega\}$
\cite{HaroL06a}. This result uses heavily that the
cocycles are over a rotation.

A system can have several dichotomies, but for the purposes of this
paper, the definition \ref{def:dichotomy} will be enough, since
we can perform the analysis presented here for each gap.

One fundamental problem for subsequent applications is the
computation of the invariant splittings (and, of course, to
ensure their existence). The computation of the invariant bundles
is closely related to the computation of iterations of the cocycle.

The first algorithm that comes to mind, is an analogue of the power method
to compute leading eigenvalues of a matrix.
Given a typical vector $(x_\theta,v)\in \E^u$, we expect that,  for
$n\gg 1$, $\M(n,\theta)v$ will be a vector in
$\E_{x_{T_\omega^n(\theta)}}^u$. Even if there are issues related 
to the $\theta$ dependence, this  may be practical if $\E^u$ is a 1-dimensional bundle.

\section{Equivalence of cocycles, reducibility}

Reducibility is a very important issue in the theory of cocycles. We have the following definition. 
\begin{defn}
We say that a cocycle $\widetilde M(\theta)$ is equivalent to
another cocycle  $M(\theta)$ if
there exists a matrix valued function
$Q :\torus^{\ell} \to GL(2d,\real)$ such that
\begin{equation}
\widetilde M(\theta) = Q(\theta +\omega)^{-1} M(\theta) Q(\theta).
\label{equivalence}
\end{equation}
\end{defn}

It is easy to check that $\widetilde M$ being equivalent to $M$ is
an equivalence relation.

If $\widetilde M$ is equivalent to a constant cocycle (i.e. independent of $\theta$), we say that
$\widetilde M$ is ``reducible.''

The important point is that, when \eqref{equivalence} holds, we have
\begin{equation}
\widetilde{\M}(n,\theta) = Q(\theta + n\omega)^{-1} \M (n,\theta)
Q(\theta). \label{iterations}
\end{equation}

In particular, if $M$ is a constant matrix, we have
\begin{equation*}
\widetilde \M (n,\theta) = Q^{-1} (\theta + n\omega) M^n Q(\theta),
\end{equation*}
so that the iterations of reducible cocycles are very easy to compute.

We will also see that one can alter the numerical stability
properties of the iterations of cocycles by choosing appropriately
$Q$. In that respect, it is also important to mention the concept of
``quasi-reducibility'' introduced by Eliasson \cite{Eliasson01}.

\section{Algorithms for fast iteration of cocycles over rotations}
\label{cocyclefast}

In its simplest form,  the algorithm for fast iteration of cocyles
is:

\begin{algorithm}[Iteration of cocycles $1$]\label{alg:doubling}
Given $M(\theta)$, compute
\begin{equation}
\widehat M(\theta) = M(\theta +\omega) M(\theta). \label{doubling}
\end{equation}
Set $\widehat M \to M$, $2\omega\to \omega$ and iterate the
procedure.
\end{algorithm}

We refer to $\widehat M$ as the renormalized cocycle and the
procedure as a renormalization procedure.

The important property is that applying $k$ times the
renormalization procedure described in Algorithm \ref{alg:doubling} amounts to compute the cocycle
$\M (2^k ,\theta)$.

Then, if we discretize the matrix $M(\theta)$ taking $N$ points (or $N$
Fourier modes) and denote by $C(N)$ the number of operations
required to perform a step of Algorithm~\ref{alg:doubling},  we
can compute $2^k$ iterates at a cost of $k C(N)$ operations.

Notice that the important point is that the cost of computing $2^k$
iterations is proportional to $k$. Of course, the proportionality
constant depends on $N$. The form of this dependence on $N$ depends
on the details on how we compute the shift and the product of matrix
valued functions.

There are several alternatives to perform the transformation
\eqref{doubling}. The main difficulty arises from the fact that, if
we have points on a equally  spaced grid, then $\theta+\omega$ will
not be in the same grid. We have at least three alternatives:
\begin{enumerate}
\item
Keep the discretization by points in a grid and  compute $M(\theta
+\omega)$ by interpolating with nearby points.
\item
Keep the discretization by points in a grid but note that the shift
is diagonal in Fourier space. Of course, the multiplication of the
matrix is diagonal in real space.
\item
Keep the discretization in Fourier space but use the Cauchy formula
for the product.
\end{enumerate}
The cost factor of each of these alternatives is, respectively,
\begin{equation}\label{costs}
\begin{split}
&C_1 (N) = O(N),\\
&C_2 (N) = O(N\log N),\\
&C_3 (N) = O(N^2).
\end{split}
\end{equation}
Besides their cost, the above algorithms
may have different stability and roundoff properties. We are not
aware of any study of these stability or round-off properties. The properties
of interpolation depend on the dimension.

In each of the cases, the main idea of the method is to precompute some blocks of
the iteration, store them and use them in future iteration. One can
clearly choose different strategies to group the blocks. Possibly, different methods can lead to different
numerical (in)stabilities. At this moment, we lack a definitive theory of
stability which would allow us to choose the blocks.


Next, we  will present an alternative consisting of using the $QR$
decomposition for the iterates. As described, for instance in
\cite{Oseledec68,EckmannR85,DieciV02},
the $QR$ algorithm seems to be rather stable to compute iterates.
One advantage is that, in the case of several gaps, it can compute
all the eigenvalues in a stable way.

\begin{algorithm}[Computation of cocycles with $QR$]\label{alg:QRdoubling}
Given $M(\theta)$ and a $QR$ decomposition of $M(\theta)$,
\[
M(\theta)=Q(\theta)R(\theta),
\]
perform the following operations:
\begin{itemize}
\item Compute $S(\theta)=R(\theta + \omega) Q(\theta)$
\item Compute pointwise a $QR$ decomposition of $S$, $S(\theta)=\bar Q(\theta) \bar R(\theta)$.
\item Compute $\widetilde Q(\theta)=Q(\theta+\omega) \bar Q (\theta)$
\item[{}] \qquad\qquad $\widetilde R(\theta)=\bar R(\theta) R (\theta+\omega)$
\item[{}] \qquad\qquad $\widetilde M(\theta)=\widetilde Q (\theta) \widetilde R (\theta)$
\item Set  $M \gets \widetilde M$
\item[{}] \quad \quad $R \gets \widetilde R$
\item[{}] \quad \quad $Q \gets \widetilde Q$
\item[{}] \quad \quad $2 \omega \gets \omega$
\end{itemize}
and iterate the procedure.
\end{algorithm}

Since the $QR$ decomposition is a fast algorithm, the total cost of
the implementation depends on the issues previously discussed (see costs in \eqref{costs}).
Instead of using $QR$ decomposition, one can also perform a $SVD$ decomposition (which is somewhat
slower).

In the case of one-dimensional maps, one can be more precise in the description of the method. Indeed, if the frequency $\omega$ has
a continued fraction expansion
$$\omega = [a_1, a_2,\ldots,a_n,\ldots],$$
it is well known
that the denominators $q_n$ of the convergents of
$\omega$ (i.e. $p_n/q_n = [a_1,\ldots, a_n]$) satisfy
\begin{equation*}
\begin{split}
&q_n = a_n q_{n-1} + q_{n-2},\\
&q_1 = a_1,\\
&q_0 =1.
\end{split}
\end{equation*}

As a consequence, we can consider the following algorithm for this particular case:
\begin{algorithm}[Iteration of cocycles $1D$]\label{contfract}
Given $\omega = [a_1,\ldots,a_n,\ldots]$ and the cocycle over
$T_\omega$ generated by $M(\theta)$, define
$$\omega^0 =0,\qquad \omega^1=\omega,\qquad
M^0 (\theta) = \Id,\qquad
M^1 (\theta)=M(\theta + (a_1-1)\omega)\cdots M(\theta).$$
Then, for $n\ge2$
$$M^{(n)} (\theta) = M^{(n-1)} (\theta + (a_n-1) \omega^{n-1})
\cdots  M^{(n-1)} (\theta) M^{(n-2)}(\theta)$$
is a cocycle over
$$\omega^{n-1} = a_n \omega^{n-1} + \omega^{n-2}$$
and we have
$$\M (q_n,\theta) = \M^{(n)} (1,\theta).$$
\end{algorithm}

The advantage of this method is that the effective rotation is
decreasing to zero so that the cocycle is becoming in some ways
similar to the iteration of a constant matrix. This method is somehow reminiscent of some algorithms that
have appeared in the mathematical literature
\cite{Rychlik92,Krikorian99a,Krikorian99b}.

\section{The ``straddle the saddle'' phenomenon and
preconditioning}\label{sec:straddlesaddle}

The iteration of cocycles has several pitfalls compared with the
iteration of matrices. The main complication from the numerical
point of view is that the (un)stable bundle does depend on the base
point.

In this section we describe a geometric phenomenon that causes some
instability in the iteration of cocycles. This instability --which
is genuine-- becomes catastrophic when we apply some of the fast
iteration methods described in Section~\ref{cocyclefast}. The
phenomenon we will discuss was already observed in \cite{HaroL06c}.

Since we have the inductive relation,
$$\M (n,\theta) = \M (n-1,\theta +\omega) M(\theta),$$
we see that we can think of computing $\M (n,\theta)$ by applying
$\M(n-1,\theta +\omega)$ to the column vectors of $M(\theta)$.

The $j^{th}$-column of $M$, which we will denote by $m_j(\theta)$, can be thought
geometrically as an embedding from $\torus^{\ell}$ to $\real^{2d}$ and is given by $M(\theta)e_j$ where $e_j$ is the $j^{th}$ vector of the canonical basis of $\real^{2d}$. If
the stable space of $\M(n-1,\theta+\omega)$ has codimension $\ell$ or
less, there could be points $\theta^* \in \torus^{\ell}$ such
that $m_j(\theta^*) \in \E_{x_{\theta^*}}^s$ and such that for every $\theta \neq \theta^*$ we have $m_j(\theta) \notin \E_{x_\theta}^s$.

Clearly, the quantity
$$\M (n-1,\theta^* + \omega) m_j(\theta^*)$$
decreases exponentially. Nevertheless, for all $\theta$ in a
neighborhood of $\theta^*$ such that $\theta \neq \theta^*$
$$\M (n-1,\theta +\omega) m_j(\theta)$$
will grow exponentially. The direction along which the
growth takes place depends on the projection of
$ m_j(\theta)$ on $\E_{x_{\theta+\omega}}^u$.

\begin{figure}[h]
 \begin{center}
 \includegraphics[width=90mm]{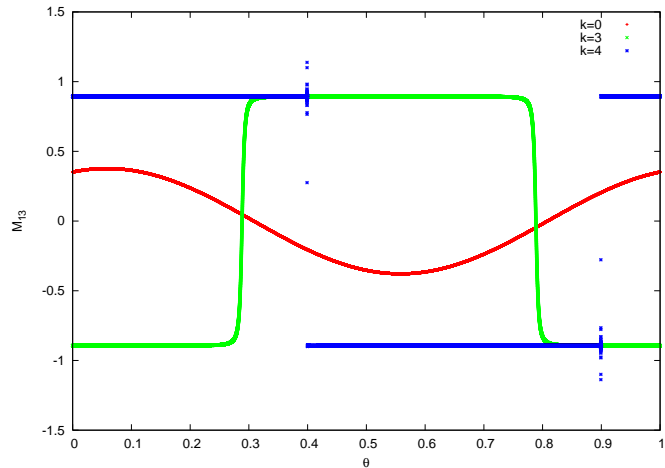}
 \caption{The straddle the saddle phenomenon. We plot one of the components of the cocycle $\M (2^k,\theta)$ for the values $k=0,3,4$. The case $k=0$ was scaled by a factor $200$.}
 \label{fig:straddlesaddle}
  \end{center}
 \end{figure}

For example, in the case $d=2$, $\ell=1$ and the stable and
unstable directions are one dimensional, the unstable components
will have different signs and the vectors $\M (n-1,\theta+\omega)
m_j(\theta)$ will align in opposite directions. An illustration of
this phenomenon happens in Figure~\ref{fig:straddlesaddle}.

The transversal intersection of the range of $m_j(\theta)$ with
$\E^s$ is indeed a true phenomenon, and it is a true
instability.

Unfortunately, if $m_j(\theta)$ is very discontinuous as a function of $\theta$, the
discretization in Fourier series or the interpolation by splines
will be extremely inaccurate so that the
Algorithm~\ref{alg:doubling} fails to be relevant.

This phenomenon is easy to detect when it happens because the
derivatives grow exponentially fast in some localized spots.

One important case where the straddle the saddle is unavoidable is
when the invariant bundles are non-orientable. This happens
near resonances (see \cite{HaroL07}). In \cite{HaroL07}, it is shown
that, by doubling the angle the case of resonances can be studied
comfortably because then, non-orientability is the only obstruction
to the triviality of the bundle.

\subsection{Eliminating the ``straddle the saddle''
 in the one-dimensional case}
Fortunately, once the phenomenon is detected, it can be eliminated.
The main idea is that one can find an equivalent cocycle which does
not have the problem (or presents it in a smaller extent).

In more geometric terms we observe that, even if the stable and
unstable bundles are geometrically natural objects, the
decomposition of a matrix into columns is coordinate dependent.
Hence, if we choose some coordinate system which is reasonably close
to the stable and unstable manifolds and we denote by $Q$ the change
of coordinates, then the cocycle
$$\widetilde {\M}(\theta)= Q(\theta +\omega)^{-1} M(\theta) Q(\theta),$$
is close to a constant. Remark that this is true only in the one-dimensional case. The picture is by far
more involved when the bundles have higher rank.

This may seem somewhat circular, but the circularity can be broken
using continuation methods. Given a cocycle which is close to
constant, fast iteration methods work and they allow us to compute
the splitting. Then if we have computed $Q$ for some $M$, we can
use it to precondition the computation of neighboring $M$.

The algorithms for the computation of bundles will be discussed
next. 


\subsection{Computation of rank-1 stable and unstable bundles using iteration of cocycles}\label{subsec:ctdyn}

The algorithms described in the previous section provide a fast way
to iterate the cocycle. We will see that this iteration method,
which is similar to the power method, gives the dominant eigenvalue $\lambda_{max}
(\theta)$ and the corresponding eigenvector $m(\theta)$.

The methods based on iteration rely strongly on the fact that the
cocycle has one dominating eigenvalue which is simple. 

Consider that we have performed $k$ iterations of the cocycle (of
course we perform scalings at each step) and we have computed
$\M(n,\theta)$, with $n=2^k$. Then, one can easily read the dominant
rank-1 bundle from the $QR$ decomposition of the cocycle
$\M (n,\theta)$, just taking the column of $Q$ associated to the
largest value in the diagonal of the upper triangular matrix $R$.
One obtains a vector $m(\theta + 2^{k} \omega)$ (and
therefore $m(\theta)$ by performing a shift of angle $-2^{k}
\omega$) of modulus 1 spanning the unstable manifold. Since,
\[ M(\theta) m (\theta) = \lambda_{max}(\theta) m (\theta + \omega),\]
we have then
\[ \lambda_{max}(\theta)=([M(\theta)m(\theta+\omega)]^T [M(\theta)m(\theta+\omega)])^{1/2}.\]

As it is standard in the power method, we perform scalings at each
step dividing all the entries in the matrix $M(\theta)$ by
the maximum value among the entries of the matrix.


Hence, for the simplest case that there is one dominant eigenvalue,
the method produces a section $m$ (spanning the unstable subbundle)
and a real function $\lambda_{max}$, which represents the dynamics on
the rank 1 unstable subbundle, such that
\[ M(\theta) m (\theta) = \lambda_{max} (\theta) m (\theta + \omega).\]

Following \cite{HaroL06b}, under certain non-resonant conditions
which are satisfied in the case of the stable and unstable
subspaces, one can reduce the 1-dimensional cocycle associated to
$M$ and $\omega$ to a constant. Hence, we look for a positive
function $p(\theta)$ and a constant $\mu$ such that
\begin{equation}\label{eq:cocycle1d}
\lambda_{max}(\theta)p(\theta)=\mu p (\theta + \omega).
\end{equation}

If $\lambda_{max}(\theta)>0$, we take logarithms on both sides of the equation
\eqref{eq:cocycle1d}. This leads to
\[ \log \lambda_{max}(\theta) + \log p(\theta)= \log \mu + \log p (\theta + \omega),\]
and taking $\log \mu$ to be the average of $\log \lambda_{max}(\theta)$ the
problem reduces to solve for $\log p(\theta)$. The case $\lambda_{max}(\theta)<0$ is analogous. Of
course, $p(\theta)$ and $\mu$ can be obtained just
exponentiating.

\section*{Acknowledgements}
The work of
R.L. and G.H. has been partially supported by NSF grants. G.H. has
also been supported by the Spanish Grant MTM2006-00478 and the
Spanish Fellowship AP2003-3411.

We thank \'A Haro, C. Sim\'o for several discussions and for
comments on the paper. The final version was written while we were
visiting CRM during the Research Programme \emph{Stability and
Instability in Mechanical Systems (SIMS08)}, for whose hospitality
we are very grateful. YS would like to thank the hospitality of the department of Mathematics of 
University of Texas at Austin, where part of this work were carried out. 
\bibliographystyle{alpha}

\bibliography{referencies}

\end{document}